\documentclass[a4paper,oneside,11pt]{amsart}
\usepackage{amssymb}
\usepackage{amsmath, amsfonts}
\usepackage{mathrsfs}
\usepackage[dvips]{graphicx}
\usepackage{psfrag}
\usepackage{hyperref}
\theoremstyle{plain}
\newtheorem{theorem}{Theorem}[section]
\newtheorem{proposition}[theorem]{Proposition}
\newtheorem{lemma}[theorem]{Lemma}
\newtheorem{remark}[theorem]{Remark}
\newtheorem{corollary}[theorem]{Corollary}

\numberwithin{equation}{section}

\newcommand{\spec}{\operatorname{{Spec}}}

\newcommand{\SA}{{\mathcal A}}

\newcommand{\hrc}{{\it hRc}}
\newcommand{\ff}{{\mathcal{F}}(f)}
\newcommand{\fg}{{\mathcal{F}}(g)}

\newcommand{\fft}{{\mathcal{F}}({f_\theta})}

\newcommand{\mm}{\medskip}
\newcommand{\no}{\noindent}
\begin{document}

%%%%%%%%%%%%%%%%%%%%%%%%%%%%%%%%%%%%%%%%%%%%%%%%%%%%%%%%%%
%%%%%%%%%%%%%%%%%%%%%%%%%%%%%%%%%%%%%%%%%%%%%%%%%%%%%%%%%

\title[Foliations and Polynomial Diffeomorphisms of $\mathbb{R}^{3}$]
      {Foliations and Polynomial Diffeomorphisms of $\mathbb{R}^{3}$}

\author{Carlos Gutierrez}\thanks{The first author was supported by FAPESP of Brazil
Grant 03/03107-9}

\author{Carlos Maquera}\thanks{The second author was supported by FAPESP of Brazil
Grant 03/03107-9}


\date{\today}

\address{Carlos Gutierrez and Carlos Maquera\\
 Universidade de S{\~a}o Paulo - S{\~a}o
Carlos \\Instituto de Ci{\^e}ncias
Matem{\'a}ticas e de Computa\c{c}{\~a}o\\
Departamento de Matem{\'a}tica\\
Av. do Trabalhador S{\~a}o-Carlense 400 \\
13560-970 S{\~a}o Carlos, SP\\
Brazil}
 \email{gutp@icmc.usp.br}
 \email{cmaquera@icmc.usp.br}
\maketitle

 \begin{abstract}
 Let $Y=(f,g,h):\mathbb{R}^{3} \to \mathbb{R}^{3}$ be a $C^{2}$ map and
 let $\spec(Y)$ denote the set of eigenvalues of the derivative
 $DY_p$, when $p$ varies in $\mathbb{R}^3$.
 We begin proving that if,
 for some $\epsilon>0,$ $\spec(Y)\cap (-\epsilon,\epsilon)=\emptyset,$
 then the foliation $\mathcal{F}(k),$   with $k\in \{f,g,h\},$
 made up by the level surfaces
 $\{k={\rm constant}\},$ consists just of planes.
 As a consequence, we prove a bijectivity result related to the
 three-dimensional case of  Jelonek's Jacobian Conjecture for
 polynomial maps of $\mathbb{R}^n.$
 \end{abstract}

 \medskip
 \medskip

 \thispagestyle{empty}
 %\tableofcontents

 %%%%%%%%%%%%%%%%%%%%%%%%%%%%%%%%%%%%%%%%%%%%%%%%%%%%%%%%%%%%%%%%%%%
 %%%%%%%%%%%%%%%%%%%%%%%%%%%%%%%%%%%%%%%%%%%%%%%%%%%%%%%%%%%%%%%%%%
 \section{Introduction}
 %%%%%%%%%%%%%%%%%%%%%%%%%%%%%%%%%%%%%%%%%%%%%%%%%%%%%%%%%%%%%%%%%%%
 %%%%%%%%%%%%%%%%%%%%%%%%%%%%%%%%%%%%%%%%%%%%%%%%%%%%%%%%%%%%%%%%%%

 Let $Y=(f,g,h):\mathbb{R}^{3} \to \mathbb{R}^{3}$ be a $C^{2}$ map and
 let $\spec(Y)$ be the set of (complex) eigenvalues of the derivative
 $DY_p$ when $p$ varies in $\mathbb{R}^3$. If
 for all $p\in \mathbb{R}^{3}$, $DY_{p}$ is non singular, (that is,
 $0\notin\spec(Y)$) then it follows from the inverse function theorem that:

 \textit{for each $k\in \{f,g,h\},$ the level surfaces $\{k={\rm constant}\}$
 make up a codimension one $C^{2}$-foliation $\mathcal{F}(k)$ on
 $\mathbb{R}^{3}.$} Our first result is the following

 \mm
 \begin{theorem}
 \label{thm:fol}
 If, for some $\epsilon>0,$ $\spec(Y)\cap (-\epsilon,\epsilon)=\emptyset,$
 then $\mathcal{F}(k),\ k\in \{f,g,h\},$ is a foliation by planes.
 Consequently, there is a foliation $F_{k}$ in $\mathbb{R}^{2}$
 such that $\mathcal{F}(k)$ is conjugate to the product of $F_{k}$ by $\mathbb{R}.$
 \end{theorem}

 To state our next results, we need to introduce some  concepts.
 Let $Y:M\rightarrow N$ be a continuous map of locally compact
 spaces. We say that the mapping $Y$ is \textit{not proper at a point} $y\in N$,
 if there is no neighborhood $U$ of the point $y$ such that the
 set $Y^{-1}(\overline{U}))$ is compact.

 The set $S_Y$ of points at which the map $Y$ is not proper
 indicates how the map $Y$ differs from a proper map. In particular
 $Y$ is proper if and only if this set is empty. Moreover, if
 $Y(M)$ is open, then $S_Y$ contains the border of the set $Y(M)$. The set
 $S_Y$ is the minimal set $S$ with a property that the mapping
 $Y:M\setminus Y^{-1}(S)\rightarrow N \setminus S$ is proper.

 Jelonek proved in \cite{Jelonek} that: if $Y:\mathbb{R}^n \rightarrow \mathbb{R}^n$ is a
 real polynomial mapping with nonzero Jacobian everywhere and
 $\textrm{codim}(S_Y) \geq 3,$ then $Y$ is a bijection (and consequently $S_Y=\emptyset$).

 On the other hand, the example of Pinchuk (see \cite{Pinchuk}) shows that there
 are real polynomial mappings, which are not injective, with nonzero Jacobian everywhere and
 with $\textrm{codim}(S_Y)=1$. Hence the only interesting case is that of
 $\textrm{codim}(S_Y)=2$ and we can state:

 \mm
 \no
 \textbf{Jelonek's Real Jacobian Conjecture.} \emph{Let
 $Y:\mathbb{R}^n \rightarrow \mathbb{R}^n $ be a real polynomial
 mapping with nonzero Jacobian everywhere. If\, ${\rm codim}(S_Y) \geq 2$
 then $Y$ is a bijection (and consequently $S_Y=\emptyset$).}

 \mm
 Jelonek \cite{Jelonek} proved that his conjecture is true in dimension two.
 Consequently, the first interesting case is $n=3$ and  $\dim (S_Y)=1$.

 \mm
 Jelonek's Real Jacobian Conjecture is closely connected with the
 following famous Keller Jacobian Conjecture:

 \mm
 \no
 \textbf{Jacobian Conjecture.}
 \textit{Let $Y:\mathbb{C}^n \rightarrow
 \mathbb{C}^n$ be a polynomial mapping with nonzero Jacobian
 everywhere, then $Y$ is an isomorphism.}

 \mm
 More precisely, Jelonek proved in \cite{Jelonek} that his Real Jacobian Conjecture
 in dimension $2n$ implies the Jacobian Conjecture in (complex) dimension
 $n$. The corresponding Jelonek's arguments and some well known
 results (\cite{BCW}, \cite{C-N}, \cite{Ess}, \cite{Yag}) will be used to
 obtain in section 3 the following version of the Reduction Theorem

 \begin{theorem}
 \label{Keller}
 Let $X_i:\mathbb{C}^n \rightarrow \mathbb{C}$
 denote the canonical i-coordinate function. If $F$, with ${\rm codim}(S_F) \geq 2$,
 is injective for all $n \geq 2$ and all polynomial maps $F:\mathbb{R}^n \rightarrow
 \mathbb{R}^n$ of the form
 $$F=(-X_1+H_1,-X_2+H_2,\ldots,-X_n+H_n)$$
 where each $H_i:\mathbb{R}^n \to \mathbb{R}$
 is either zero or homogeneous of degree 3, and the Jacobian matrix
 $JH$ (with $H=(H_1,H_2,\ldots,H_n)$) is nilpotent, then the Jacobian Conjecture is true.
 \end{theorem}

 Notice that in theorem above $\spec (F)=\{-1\}$.

 \mm

 Related with Theorem \ref{Keller} and Jelonek conjecture we prove
 the following.

 \begin{theorem}
 \label{theo3Jelonek}
 Let $Y=(f,g,h):\mathbb{R}^3 \rightarrow \mathbb{R}^3$ be a polynomial map
 such that $\spec(Y) \cap [0,\varepsilon)=\emptyset$, for some
 $\varepsilon>0$. If ${\rm codim}(S_Y) \geq 2$ then $Y$ is a bijection.
 \end{theorem}

 This result partially extends also the bi-dimensional results of \cite{CGL}
 and \cite{agr} (see also \cite{Cam} -- \cite {Cha}, \cite{fe-gu},
 \cite{Gu} -- \cite{N-X}, \cite{CGPV}).
%\cite{CGPV}, \cite{G-R}, \cite{GV1}, \cite{GV2}, \cite{GV3}

\section{Half-Reeb Components and The Spectral Condition}

 Let us recall the definition of a vanishing cycle stated in
 conformity  with our needs. Let $Y=(f,g,h):\mathbb{R}^3 \to \mathbb{R}^3$
 be a $C^{2}$ map such that, for all $p\in \mathbb{R}^3$, $DY_p$ is non-singular.
 Given $k\in\{f,g,h\},$ a \textit{vanishing cycle}
 for the foliation $\mathcal{F}(k)$ is a $C^2$-embedding
 $f_{0}:S^{1}\to \mathbb{R}^3$ such that:
 \begin{itemize}
    \item[(a)] $f_{0}(S^{1})$ is contained in  a leaf $L_{0}$ but it
    is not homotopic to a point in $L_{0}$;
    \item[(b)] $f_{0}$ can be extended to a $C^2-$embedding
    $f:[0,1]\times S^{1}\to \mathbb{R}^3,\ f(t,x)=f_{t}(x),$ such that for all
    $t>0$, there is a 2-disc $D_t$ is contained in a leaf $L_{t}$,
    such that $\partial D_t=f_t(S^1)$;
    \item[(c)] for all $x\in S^{1},$ the curve $t\mapsto f(t,x)$
    is transversal to the foliation $\mathcal{F}(k)$ and, for all $t
    \in (0,1),\,D_t$ depends continuously on $t$.
 \end{itemize}

 We say that the leaf $L_{0}$ \textit{supports} the vanishing cycle
 $f_{0}$ and that $f$ is the map associated to $f_0$.

 The \textit{half-Reeb component for \/} $\mathcal{F}(k)$
 (or simply the $\hrc$ for $\mathcal{F}(k)$)
 associated to the
 vanishing cycle $f_0$ is the region
 $$
 \mathcal{A}=\left(\bigcup_{t \in (0,1]} \, D_t \right)\cup L \cup
 f_0(S^1)
 $$
 where $L$ is the connected component of $L_0-f_0(S^1)$ contained in
 the closure of $\cup_{t \in (0,1]}\,D_t$. The transversal section $A=f([0,1]{\times}S^1)$
 to the foliation $\mathcal{F}(k)$
 is called the \textit{compact face \/} of $\mathcal{A}$ and the leaf $L \cup
 f_0(S^1)$ of $\mathcal{F}(k)|_{\mathcal{A}}$ is called the
 \textit{non-compact face \/} of $\mathcal{A}$.

 \begin{figure}[ht!]
 \label{figi:halfrc}
 \psfrag{L}{\footnotesize $L$}
 \psfrag{A}{\footnotesize$A$}
 \psfrag{g}{\footnotesize$f_0(S^1)$}
       \includegraphics[width=4.5cm,height=4.5cm,keepaspectratio]{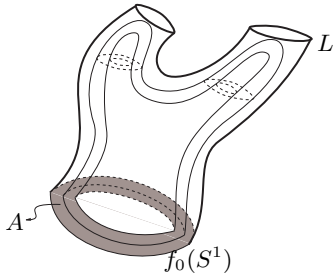}
 \caption{A half-Reeb component.}
 \end{figure}

 \mm
 \begin{remark}
 \noindent
 \begin{enumerate}
    \item {\rm It will be seen in Proposition \ref{prop:cevanes}, that if
    $\mathcal{F}(k)\,,\ k\in\{f,g,h\},$ has a leaf which is not homeomorphic to the plane,
    then $\mathcal{F}(k)$ has a half-Reeb component. }
    \item {\rm The connection between half-Reeb components and the spectral condition
    on $Y$ (that is, $\spec(Y)\cap (-\varepsilon,\varepsilon)=\emptyset$)
    is given by Theorem \ref{thm:fol}.}
 \end{enumerate}

 \end{remark}

 \mm
 The following proposition is obtained by using classical arguments of Foliation
 Theory (see \cite{Camacho} and \cite{godbillon}). For sake of completeness we give the main lines of its
 proof. Let $D^2=\{(x,y)\in\mathbb{R}^2:x^2+y^2\leq1\}$ denote the closed
 2-disc.

 \mm
 \begin{proposition}
 \label{prop:cevanes}
 If $\mathcal{F}(k), \,$ with $ k\in\{f,g,h\},$ has a leaf $L$ which is not
 homeomorphic to the plane, then $\mathcal{F}(k)$ has a vanishing
 cycle.
 \end{proposition}
 \begin{proof}
 Let $\eta:S^{1}\to L$ be an embedding which is
 not null homotopic in $L.$ Since $\eta$ is null homotopic in $\mathbb{R}^{3},$
 we may extend it to a $C^{2}$-immersion $\eta:D^2 \to \mathbb{R}^3$,
 which is in general position with respect to $\mathcal{F}(k)$. It
 this way we are supposing that the contact set $C_\eta$, made up by
 the points of $D^2$ at which $\eta$ meets tangentially
 $\mathcal{F}(k)$, is finite and is contained in $D^2\setminus S^1$.

  Via $\eta$, the foliation $\mathcal{F}(k)$ induces a foliation $\mathcal{G}$
 (with singularities) on $D^{2}$. We claim that it is possible construct a
  vector field $G$ on $D^{2}$ such that the foliation $\mathcal{G}$
  is induced by $G.$ In fact, as $\eta$ is in general position with respect to
  $\mathcal{F}(k)$, the foliation $\mathcal{G}$ has finitely many singularities
  each of which is locally topologically equivalent either to a center or to a
  saddle point of a vector field. This implies that $\mathcal{G}$ is locally
  orientable everywhere. As $D^2$ is simple connected, $\mathcal{G}$
  is globally orientable. This proves the existence of the vector
  field $G$. Certainly, we may assume that $\eta$ has been chosen so
  that no pair of singularities of
  G is taken by $\eta$ into the same leaf of $\mathcal{F}(k);$
  in other words, $G$ has no saddle connections.

  We claim that $G$ has no limit cycles. In fact, otherwise,
  the Poincar{\'e}-Bendixon theorem would imply that there is a orbit of
  $G$ which spirals towards a limit cycle $C.$
  Hence, the leaf of $\mathcal{F}(k)$ containing $C$ would have a
  non trivial holonomy group. This contradiction proves our claim.

  \begin{figure}[!ht]
  \label{fig:evanescente}
   \begin{center}
      \psfrag{a}{\footnotesize ($a$)}
      \psfrag{b}{\footnotesize ($b$)}
      \psfrag{c}{\footnotesize ($c$)}
      \psfrag{D}{\footnotesize $D$}
      \psfrag{g}{\footnotesize $\gamma_i$}
      \psfrag{s}{\footnotesize $s_i$}
      \psfrag{c2}{\footnotesize $c_i$}
      \psfrag{Di}{\footnotesize $D_i$}
              \includegraphics[width=13cm,height=13cm,keepaspectratio]{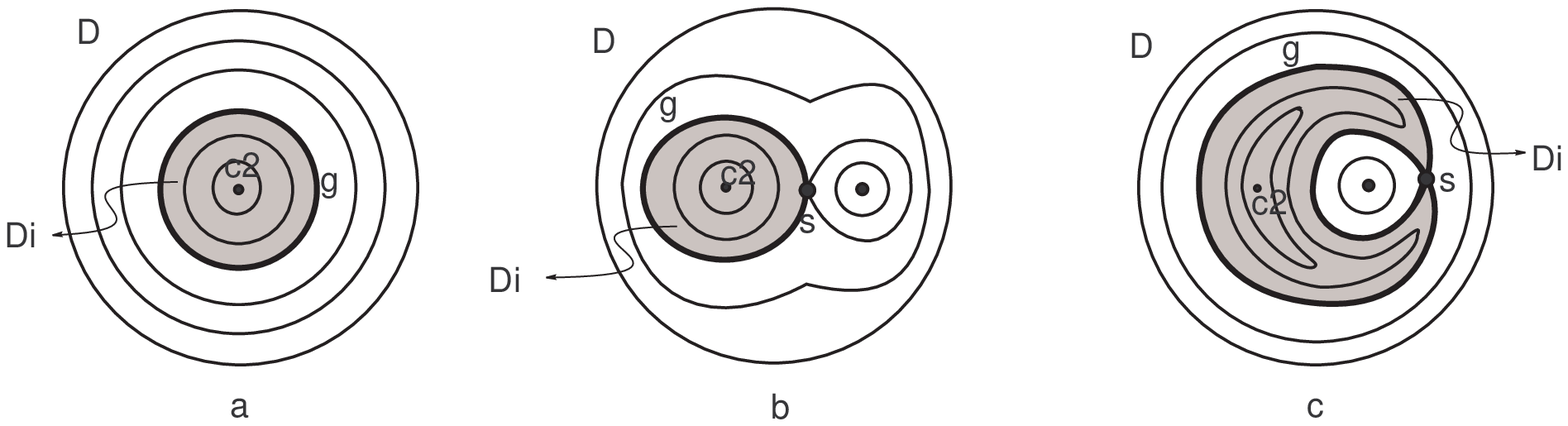}
              \caption{\footnotesize }
    \end{center}
 \end{figure}

 Let $c_{1},\dots,c_{\ell}$ be the center singularities of $G$.
 Given $i\in\{1,\dots,\ell\}$, there exists a G-invariant open
 2-disc $D_{i}\subset D^2$ such that:

 \begin{itemize}
    \item[(a1)] $c_i \in D_i$ and every orbit of $G$ passing through
    a point in $D_{i}\setminus \{c_{i}\}$ is a closed orbit;

    \item[(a2)] for every closed orbit $\gamma\subset D_{i}$ of $G$,
    $\eta(\gamma)$ is homotopic to a point in its corresponding leaf
    of $\mathcal{F}(k)$.

    \item[(a3)] the 2-disc $Di$ is the biggest one satisfying
    properties (a1) and (a2) above.
 \end{itemize}

 \medskip

Notice that the frontier $\gamma_i$ of $D_i$ has to be G-invariant.
We claim that

\medskip

\begin{itemize}
\item[(b)] If, for some $i \in \{1,2,\ldots,\ell\}$, $\gamma_i$ is a
    closed orbit of G, then $\eta(\gamma_i)$ is a vanishing cycle (see Fig.
    \ref{fig:evanescente}(a)), and the proposition is proved.
 \end{itemize}

In fact, if $\gamma_i$ is a closed orbit of G such that
$\eta(\gamma_i)$ is homotopic to a point in its corresponding leaf,
then, by a well known result of foliation theory, there exists a
neighborhood $V_i\subset D^2$ of $\gamma_i$ such that the image by
$\eta$ of every orbit of $G$, contained in $V_i$, is homotopic to a
point in its corresponding leaf. This contradiction with the
maximality of $D_i$ proves (b).

\medskip

Therefore, we may suppose, from now on, that:

\medskip

\begin{itemize}
    \item[(c)] for every $i\in \{1,\dots,\ell\},$ $\gamma_i$ is either the union
    of a saddle singularity $s_{i}$ of $G$ and one of its separatrices or the union
    of a saddle singularity $s_{i}$ and its two separatrices, see (b) and (c) of
    Figure \ref{fig:evanescente}.
\end{itemize}

\medskip

By studying the phase portrait of $G$, we may conclude that

\medskip

\begin{itemize}
    \item[(d)] if (b) is not satisfied, there must exist $i \in
    \{1,2,\ldots,\ell\}$ such that $\gamma_i$ is the union of a
    saddle singularity $s_i$ of $G$ and one of its separatrices.
\end{itemize}

\mm

We claim that:

\mm

\begin{itemize}
    \item[(e.1)] If $\eta(\gamma_i)$ is homotopic to a point in its
    corresponding leaf, then $\eta$ can be deformed to a
    $C^{2}$-immersion $\tilde{\eta}:D^2 \rightarrow \mathbb{R}^3$
    which is in general position with respect to $\mathcal{F}(k)$ and
    such that $\# C_{\tilde{\eta}} < \# C_{\eta}$;
\end{itemize}

\mm

\begin{itemize}
    \item[(e.2)] If $\eta(\gamma_i)$ is a not homotopic to a point
    in its corresponding leaf, then $\eta$ can be deformed to a
    $C^{2}$-immersion $\tilde{\eta}:D^2 \rightarrow
    \mathbb{R}^3$ for which (b) above is satisfied.
\end{itemize}

\medskip

 In fact, let us prove (e.1). By using Rosenberg's arguments (see
 \cite[pag. 137]{Rosenberg}), via a deformation of $\eta$, supported
 in a neighborhood of $\overline{D}_i$, we can eliminate the saddle
 singularity $s_i$ and the center singularity $c_i$. the proof of
 (e.2) is similar and will be omitted.

 Using (e.1) as many times as necessary, it follows from (d) that we
 will arrive to the situation considered in (e.2). this proves the
 proposition.
 \end{proof}

%\mm
% \begin{remark}
% {\rm In the proof of proposition above, by deforming $\eta$, if necessary,
%  we may assume that the foliation $G$ on $D^2$ is as the Figure \ref{fig:evanescente} (a).}
% \end{remark}

 \mm
  \begin{lemma}
  \label{lem:lr}
  Let $\mathcal{F}_{i}\,,\ i=1,2,3,$ be a $C^{2}$ foliation on $\mathbb{R}^{3}$
  without holonomy such that for $j\neq i,\ \mathcal{F}_{j}$ is transversal to
  $\mathcal{F}_{i}\,.$ Let $L$ be a leaf of $\mathcal{F}_{1}\,.$ If $\mathcal{F}_{2}
  |_{L}$ is the foliation on $L$ that is induced by $\mathcal{F}_{2}\,,$ then every
  leaf of $\mathcal{F}_{2}|_{L}$ is homeomorphic to $\mathbb{R}.$
  \end{lemma}

  \begin{proof}
  Suppose that there exists a leaf $S$ of $\mathcal{F}_{2}|_{L}$, homeomorphic to
  $S^{1}.$ The fact that $\mathcal{F}_{2}$ is without holonomy and $\mathcal{F}_{3}
  |_{L}$ is transversal to $\mathcal{F}_{2}|_{L}$ implies that there exists a
  neighborhood $C$ of $S$ in $L$ such that every leaf of $\mathcal{F}_{2}|_{L}$ passing
  through a point in $C$ is homeomorphic to $S^{1}$ and is not homotopic to a point
  in $L.$ Moreover, the leaves of $\mathcal{F}_{3}|_{L}$ restricted to $C$ are curves
  starting at one connected components of $\partial C$, and ending at the other one.

  Let $D$ be a smoothly immersed open 2-disc containing $S,$ which we may
  assume to be in general position with respect to $\mathcal{F}_{3}\,.$
  Let $\mathcal{G}_{3}$ be the foliation (with singularities) of $D$ which is
  induced by $\mathcal{F}_{3}\,.$ Then, $\mathcal{G}_{3}$ is
  transversal to $S$.

  We claim that $\mathcal{G}_{3}$ has no limit cycles, otherwise,
  the Poincar\'e-Bendixon theorem implies that there is a leaf of
  $\mathcal{G}_{3}$ which spirals towards a limit cycle $\gamma.$
  Hence, the leaf of $\mathcal{F}_{3}$ containing $\gamma$ would have a
  non trivial holonomy group. This contradiction proves our claim.
  It follows from the claim above that $\mathcal{G}_{3}$, has
  exactly one singularity. Since $\mathcal{G}_{3}$ is transversal to $S$,
  this singularity is an attractor. But $D$ in general position with respect to
  $\mathcal{F}_{3}$ means that $\mathcal{G}_{3}$ has a finite number of singularities,
  each of which is either a center or a saddle point. This contradiction concludes the
  proof.
  \end{proof}

  \mm
  \begin{remark}
  {\rm
 Let $k\in\{f,g,h\}.$ As $k$ is a submersion, the foliation $\mathcal{F}(k)$  is
  without holonomy.
  }
  \end{remark}

 \mm
 \begin{corollary}
 \label{cor:mono}
  Let $\{i,j,k\}$ be an arbitrary permutation of
  $\{f,g,h\}.$ If $L$ is a leaf of $\, \mathcal{F}(i) \,$  and $l$ is a leaf
  of $\, \mathcal{F}(j)|_{L}\, $ then $k\vert_l$ is regular; in this way
  $\mathcal{F}(j)|_{L}$  and $\mathcal{F}(k)|_{L}$ are transversal to each other.
  \end{corollary}

 \mm
 For each $\:\theta\in\mathbb{R}\:$ let $T_\theta,\ S_{\theta}:\mathbb{R}^{3}\to \mathbb{R}^{3}$
 be the linear transformations defined by the matrices

 \[
 \left( \begin{array}{ccc}
          \cos \theta &  -\sin \theta & 0\\
          \sin \theta &  \phantom{+} \cos \theta & 0\\
          0& 0 & 1
        \end{array} \right)\,,
 \ \ \ \
 \left( \begin{array}{ccc}
          \cos \theta &  0 & -\sin\theta  \\
           0 &  1 & 0\\
          \sin \theta & 0 & \phantom{+} \cos \theta
        \end{array} \right),
 \]
 \mm
 respectively. Note that $T_\theta$ (resp. $S_{\theta}$) restricted
 to the $xy$-plane  (resp. $xz$-plane) is the rotation
 \[
 \left( \begin{array}{cc}
          \cos \theta &  -\sin \theta \\
          \sin \theta &  \phantom{+} \cos \theta
  \end{array} \right)
 \]

 \mm
 Let $\Pi:\mathbb{R}^3\to\mathbb{R}$ given by $\Pi(x,y,z)=x$.
 The following proposition will be needed.

 \mm
 \begin{proposition}
 \label{prop:x-projection}
 Let $Y=(f,g,h)\colon\mathbb{R}^3 \to \mathbb{R}^3$ be a $C^{2}$
 map such that $0\notin{\rm Spec}(Y)$ and $\mathcal{A}$ be a \hrc\
 of $\ff$. If \,$\Pi(\SA)$ is
 bounded, then there is an $\epsilon>0$ and
 $K_{\theta}\in\{S_{\theta},T_{\theta}\}$ such that, for all
 $\theta\in (-\epsilon,0) \cup (0,\epsilon)$,
 $\mathcal{F}({f_\theta})$ has a \hrc\ $\mathcal{A}_\theta$ such
 that $\Pi(\mathcal{A}_\theta)$ is an interval of infinite length,
 where $(f_\theta, g_\theta, h_{\theta}) = K_\theta\circ
 Y\circ K_{-\theta}\,$.
 \end{proposition}

 \begin{proof}
 If $\Pi(\SA)$ is bounded, then either $\{y:(x,y,z)\in \SA\}$
 or $\{z:(x,y,z)\in \SA\}$ is an interval of infinite length.
 We are going to show that, if $\{y:(x,y,z)\in \SA\}$
 is an interval of infinite length, then, for $K_{\theta}=T_{\theta}$,
 $\Pi(\mathcal{A}_\theta)$ is an interval of infinite length.
 The proof of the other case is analogous in which case,
 the proposition is satisfied for $K_{\theta}=S_{\theta}$. Then, assume that
 $\{y:(x,y,z)\in \SA\}$ is an interval of infinite length.

\noindent
\begin{enumerate}
\item[(a)]Let $\theta\in\mathbb{R}$ be such that, for all
$m\in\mathbb{Z}$,
    $\theta\neq\frac{m\pi}{2}$. Then $\fft$ is transversal
    to both ${T_\theta}(\ff)$ and ${T_\theta}(\fg).$
\end{enumerate}

 In fact, assume by contradiction that there exist $p\in \mathbb{R}^{3}$ such
 that $L_{T_{\theta}(p)}(f_{\theta})$ and $T_{\theta}(L_{p}(f))$ (the leaves through
 $T_{\theta}(p)$ of $\fft$ and $T_{\theta}(\ff),$ respectively)
 are tangent at $T_{\theta}(p).$ This implies that every
 $C^{1}$ curve in $T_{\theta}(L_{p}(f))$ passing through $T_{\theta}(p)$
 is tangent to $L_{T_{\theta}(p)}(f_{\theta})$ at $T_{\theta}(p)$.
 But, we will exhibit a $C^{1}$ curve $\alpha_{\theta}:(-1,1)\to T_{\theta}(L_{p}(f))$
 passing through $T_{\theta}(p)$ which is not tangent to $L_{T_{\theta}(p)}(f_{\theta})$
 at $T_{\theta}(p).$ Indeed, we consider $\alpha_{\theta}:(-1,1) \to T_{\theta}(L_{p}(f))$
 defined by $\alpha_{\theta}=T_\theta \circ \alpha$ where $\alpha:(-1,1)\to \mathbb{R}^3$ is
 a $C^{1}$ curve contained in $L_{p}(f)\cap L_{p}(h)$ with $\alpha(0)=p$ and
 $\alpha'(0)\neq 0$. By Corollary \ref{cor:mono}, $(g \circ \alpha)'(0)\ne
 0.$ Hence, as $f(\alpha(t))\equiv \mbox{constant},$
 $\sin\theta \neq 0$ and
 $$
 (f_\theta \circ \alpha_\theta)(t)=
  (\cos\theta)f(\alpha(t))-(\sin\theta)g(\alpha(t)), \ \ t\in(-1,1),
 $$
 we obtain that

 $$(f_\theta \circ \alpha_\theta)'(0)= - \sin \theta\, (g\circ\alpha)'(0) \neq 0$$

  \no and so $\alpha_{\theta}$ is not tangent to $L_{T_{\theta}(p)}(f_{\theta})$
  at $T_\theta(p)=\alpha_{\theta}(0)$. This contradiction proves that $\fft$
  is transversal to ${T_\theta}(\ff).$ Similarly we prove that $\fft$ is transversal
  to ${T_\theta}(\fg).$

 \mm
 Take $\Sigma$ diffeomorphic to the open annulus $\{z \in \mathbb{C}:1<|z|<2\}$,
 transversal to $\mathcal{F}(f)$ and containing the compact face $A$ of $\mathcal{A}$.
 Since, for $\theta$ enough small, $Y_{\theta}$ and $T_{\theta}$ are $C^{1}$ close to
 $Y$ and to the identity $T_O$, respectively, we can take $\Sigma$ so that
 \begin{itemize}
    \item[(b)] there exist $\varepsilon>0$ such that, for all $\,\theta \in
    (-\varepsilon,\varepsilon),\,T_{\theta}(\Sigma)\,$ is transversal to both
    $T_{\theta}(\ff)$ and $\fft.$
 \end{itemize}

\mm

 Let $\mathcal{G}_{\theta}$ be the foliation
 in $T_{\theta}(\Sigma)$ which is induced by $\fft$. As $\fft$ is without holonomy,
 we can take $\varepsilon > 0$ so that
 \begin{itemize}
    \item[(c)] for all $\theta \in (-\varepsilon,\varepsilon)$, there exist open
 cylinders $A_{\theta}^-,\,A_{\theta}^+\subset T_{\theta}(A_{0})$
 made up by closed trajectories of $\,\mathcal{G}_{\theta}\,$ such
 that $\,A_{\theta}^- \subset
 T_{\theta}(A),\;A^+_{\theta}\,\cap\,T_{\theta}(A)=\emptyset, \;$
 $A_{\theta}^-\cap T_{\theta}(\partial A)=\emptyset= A_{\theta}^+\cap
 T_{\theta}(\partial A) \;$ and both $\; A_{\theta}^-\;$ and
 $A_{\theta}^+\;$ are the biggest cylinders with these properties.
 \end{itemize}
 We claim that:
 \begin{itemize}
    \item[(d)] every leaf of $\mathcal{G}_{\theta}$ contained in $A_{\theta}^+$
    is not homotopic to a point in its corresponding leaf of $\fft$.
 \end{itemize}

 In fact, assume by contradiction that there exist a leaf $\gamma$ of
 $\mathcal{G}_{\theta}$ contained in $A_{\theta}^+$  and bounding
 a closed 2-disc $D(\gamma)$ contained
 in a leaf of $\fft.$
 If $L$ is the non-compact face of $\mathcal{A}$
 and $\tilde D(\gamma)\subset T_\theta(A_0)$ is the disc bounded by
 $\gamma,$ then the 2--sphere $D(\gamma)\cup (\tilde D(\gamma))$
 meets $T_\theta(L)$ at a circle contained in $\tilde D(\gamma).$
 Therefore, as the referred 2--sphere separates $\mathbb{R}^3$,
 $T_{\theta}(L)$ has to meet $D(\gamma)$ and so there exists a closed
 2-disc $D_{0}(\gamma)\subset D(\gamma)$ such that
 $\partial D_{0}(\gamma)=D(\gamma)\cap T_{\theta}(L)$. Consequently,
 there is at least one point in $D_{0}(\gamma)$
 where $T_{\theta}(\ff)$ and $\fft$ are tangent, contradicting
 (a). This proves (d).

 By using a similar argument we may also obtain that
 \begin{itemize}
    \item[(e)] every leaf of $\mathcal{G}_{\theta}$ contained in $A_{\theta}^-$
    is homotopic to a point in its corresponding leaf of $\fft$.
 \end{itemize}

 \vspace{.3cm}
 In what follows of this proof, every time that we refer to Lemma \ref{lem:lr},
 we will be assuming that it is been applied to the three foliations $\fft,$
 ${T_\theta}(\ff)$ and ${T_\theta}(\fg).$

 From (c), (e) and Lemma \ref{lem:lr}, we obtain that there exists a leaf $\gamma$
 of $\mathcal{G}_\theta$ contained in $T_{\theta}(A_{0})\setminus (A_{\theta}^{-}\cup
 A_{\theta}^{+})$ which is a vanishing cycle of $\fft$ and such that
 \begin{itemize}
    \item[(f)] $\gamma \cap T_{\theta}(\partial A)$ is a nonempty finite set.
 \end{itemize}

 \vspace{.3cm}
 Let $\mathcal{A}_\theta$ be the \hrc\ of $\fft$ with non-compact face $L_{\theta}$
 and compact face contained in $T_{\theta}(\Sigma)$ and bounded by $\gamma$. Notice
 that $L=L_O$. Let $a_{1},\dots,a_{2\ell}\in A \cap L$ be such that
 $$
 \gamma \cap T_{\theta}(A \cap L)= \{T_{\theta}(a_{1}),\dots,
 T_{\theta}(a_{2\ell})\}
 $$
 Up to small deformation of $\Sigma,$ if necessary, we may assume
 that,
 for all $i=1,\dots,2\ell,$ the connected component $\Gamma_{i}$ of
 $T_{\theta}(L)\cap L_{\theta}$ that contains $T_{\theta}(a_{i})$
 is a regular curve (not reduced to a single point).

 We claim that

 \begin{itemize}
    \item[(g)] There exists $i_0\in \{1,\dots,2\ell\}$ such that $\Gamma_{i_0}$ is
    non-compact.
 \end{itemize}

 In fact, suppose by contradiction that $\Gamma_i$ is compact for
 every $i\in \{1,\dots,2\ell\}$. Recall that $L_{\theta}$ (resp.
 $T_{\theta}(L)$) is the noncompact face of $\mathcal{A}_{\theta}$
 (resp. of $T_{\theta}(\mathcal{A})$).

 Let $U(L_{\theta})$ (resp. $U(T_{\theta}(L)$) be the unbounded
 connected component of $L_{\theta}\setminus (\cup_{i=1}^{2\ell}\Gamma_i)$
 (resp. of $T_{\theta}(L)\setminus (\cup_{i=1}^{2\ell}\Gamma_i)$). As
 $U(L_{\theta})\cap U(T_{\theta}(L))=\emptyset$ and both
 $\partial(\mathcal{A}_\theta)$ and $\partial(T_\theta(\mathcal{A}))$
 separate $\mathbb{R}^3$, we have that either
 $$
 U(L_{\theta})\subset T_{\theta}(\mathcal{A})\ \
 {\rm or}\;\,U(T_{\theta}(L)) \subset \mathcal{A}_{\theta}\,,
 $$
 respectively.
 If $U(L_{\theta})\subset T_{\theta}(\mathcal{A})$ then, since all leaves of
 $T_{\theta}(\mathcal{F}(f))|_{T_{\theta}(\mathcal{A})}$ passing through points
 in the interior of $T_{\theta}(\mathcal{A})$ are closed $2$-discs, it follows
 that $T_{\theta}(\mathcal{F}(f))|_{L_{\theta}}$ has infinitely many leaves which
 are homeomorphic to $S^1$, contradicting Lemma \ref{lem:lr}.
 Analogously, if $U(T_{\theta}(L)) \subset \mathcal{A}_{\theta}$
 we obtain a contradiction with Lemma \ref{lem:lr}.
 This proves (g).

 Now we claim that

 \begin{itemize}
 \item[(h)] $\Pi(\mathcal{A}_\theta)$ is an interval of infinite length.
 \end{itemize}

 Let $\Gamma_i=\Gamma_{i_0}$ be as in (g).
 As $\Pi(\Gamma_i) \, \subset \, \Pi(\mathcal{A}_\theta)\cap
 \Pi(T_\theta(\mathcal{A})),$ it is enough to
 prove that  $\Pi(\Gamma_i)$
 is an interval of infinite length.
 Since $\Gamma_i \, \subset \, L_\theta \, \cap \, T_\theta(L)$,
 we have that $T_{\theta}^{-1}(\Gamma_i)\subset L \subset \mathcal{A},$
 consequently $\Pi(T_{\theta}^{-1}(\Gamma_i))\subset \Pi(\mathcal{A}).$
 Now, if $\Pi(\Gamma_i)$ was bounded, then the subinterval $\Pi(T_{\theta}^{-1}
 (\Gamma_i))$ of $\Pi(\mathcal{A})$
 would have infinite length, contradicting the
 assumption
 that $\Pi(\mathcal{A})$ is bounded. This proves (h) and concludes
 the proof of this proposition.
\end{proof}

 \mm

 \begin{proof}[Proof of Theorem \ref{thm:fol}]
 By Palmeira's theorem, see \cite{pal}, it is sufficient to show that
 $\mathcal{F}(k),\ k\in\{f,g,h\},$ is a foliation by planes. Suppose
 by contradiction that $\ff$ has a leaf which is not homeomorphic to
 $\mathbb{R}^2.$ It follows, from Proposition \ref{prop:cevanes},
 that $\ff$ has a half-Reeb component $\SA$. Hereafter we will use
 the fact that existence of a half-Reeb component and the assumptions
 of Theorem \ref{thm:fol} are open in the Whitney $C^2$ topology, in
 particular we shall assume, from now on, that $Y$ is smooth. Let
 $\Pi:\mathbb{R}^3\to \mathbb{R}$ be the orthogonal projection onto
 the first coordinate. By composing with a transformation
 $T_{\theta}$ if necessary (see Proposition~\ref{prop:x-projection})
 we may assume that $\Pi(\SA)$ is an unbounded interval. To simplify
 matters, let us suppose that $[b,\infty)\subset \Pi(\SA)$ and that
 $\Pi(A) \cap [b,\infty)=\emptyset$, where $A$ is the compact face of
 $\SA$.

 By Thom's Transversality Theorem for jets \cite{golu}, we may assume
 that $\ff$ has generic contact with the foliation
 $\mathcal{F}(\Pi).$ In this way, as $f$ is a submersion,
 \begin{enumerate}
    \item[(a1)] the  contact manifold
    $T=\{(x,y,z)\in \mathbb{R}^3; f_{y}(x,y,z)=0=f_{z}(x,y,z)\}$
    is a subset of $\{(x,y,z) \in \mathbb{R}^3:f_x(x,y,z) \neq 0\}$ made up of regular
    curves;
    \item[(a2)] there is a discrete subset $\Delta$ of $T$ such that if $p\in T\setminus
    \Delta$, then $\Pi,$ restricted to the leaf of $\ff$ passing
    through $p,$ has a Morse-type singularity at $p$ which is either
    a saddle point or an extremal (maximum or minimum) point. see Figure \ref{fig:contact}.
    \end{enumerate}

 \begin{figure}[ht!]
 \label{fig:contact}
 \begin{center}
       \psfrag{CC}{\footnotesize $p\in T\setminus\Delta$: center contact}
       \psfrag{SC}{\footnotesize $p\in \Delta$: saddle contact}
       \psfrag{p}{\footnotesize $p$}
   \includegraphics[width=8.5cm,height=8.5cm,keepaspectratio]{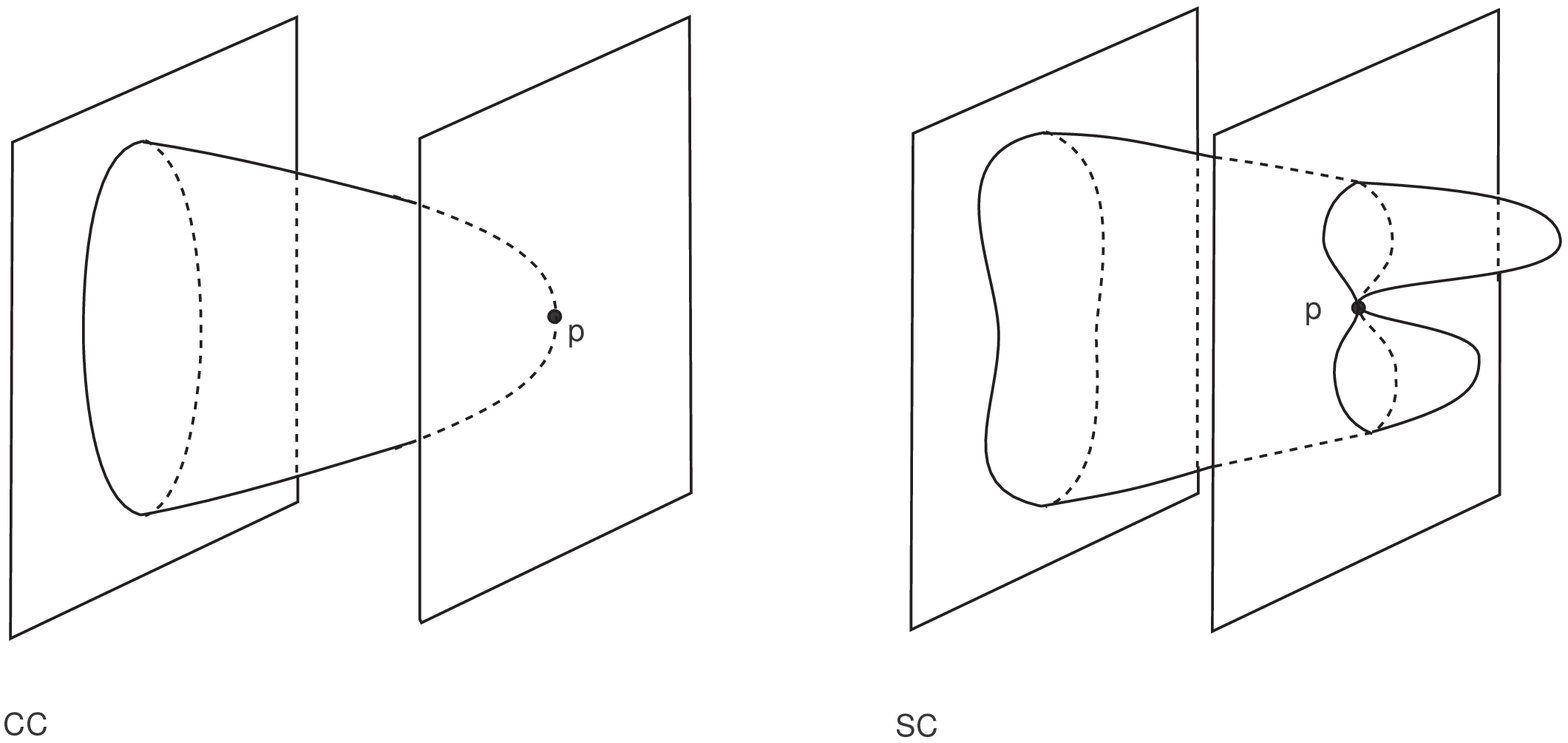}
 \caption{\footnotesize }
 \end{center}
 \end{figure}

 Then, if $a> b$ is large enough,

 \begin{enumerate}\label{maxima_dobra}
    \item[(b)]for any $x\ge a,$  the plane $\Pi^{-1}(x)$
        intersects exactly one leave $L_x \subset \SA \, $ of
        $\, \ff|_{\SA} \,$ such that $\Pi(L_x)\cap (x,\infty) =
        \emptyset.$ In other words, $x$ is the supremum of the set
        $\Pi(L_x).$ Notice that $L_x$ is a disc whose boundary
        is  contained in the compact face of $\SA$.
 \item[(c)] if $x\geq a$ then $T_x = L_{x}\cap \Pi^{-1}(x)$ is
 contained in
 $T\cap \mathcal{A}.$
   \item[(d)] if $p\in T_x$ then $p\in T\setminus \Delta$
        is a maximum point for the restriction
        $\Pi|_{L_x}\,$.
 \end{enumerate}

 Notice that $T_x$ is a finite set disjoint of $\Delta$, for every $x\geq a$.
 Hence, the map $x \in [a,\infty)\,\longmapsto\,\# T_x$ is upper semi
 continuous, were $\# T_x$ denotes the cardinal number of $T_x$.
 To motivate what is claimed in (e) below, we observe that if, for some
 $x_0\in [b,\infty)$ and for some $p\in T_{x_0}$, we had that  $\#(T_{x_0})>1$
 and $0<f_x(p)< \min \{f_x(q) : q\in T_{x_0}\setminus\{p\} \},$ then, we
 would obtain that, for some $\epsilon>0$ and for every $x\in (x_0 - \epsilon, x_0)
 \cup (x_0, x_0 + \epsilon),\, \,$ $\#T_x = 1$; in this way, there
 would exist a smooth curve $\eta:(x_0-\varepsilon,x_0+\varepsilon)\mapsto T$
 such that $\eta(x_0)=p \in T_{x_0}$ and, for all $x \neq x_0$,
 $T_x=\{\eta(x)\}$.

 Therefore, by $(b)-(d)$ and by using Thom's Transversality Theorem
 for jets, we may assume the following stronger statement:

 \begin{enumerate}
 \item[(e)] there is an increasing sequence
 $\, F=\{a_i\}_{i\ge 1}\, $
 in $\, [a,+\infty),\,$ at most countable, such that if $x\in [a,+\infty)\setminus F$,
 then  $T_{x}$ is a one-point set.
 \end{enumerate}

 If $x\in [a,+\infty)\setminus F$ and $T_x=\{(x,\eta_1(x),\eta_2(x))\}$, define
 $\eta:[a,+\infty)\setminus F\to T $ by $\eta(x)=(x,\eta_1(x),\eta_2(x)).$
 Observe that $\eta$ is a smooth embedding and, since $f|_{\mathcal{A}}$
 is continuous and bounded,

 \begin{enumerate}
 \item[(f)] $f\circ \eta$ extends continuously to a  strictly
 monotone bounded map defined in $[a,+\infty)$ such that, for
 all $x\in [a,+\infty)\setminus F$, $f_x(\eta(x))$ has constant sign.
 \end{enumerate}

 Therefore, there exists a real constant $K$ such that
 $$
 \begin{array}{ccl}
   K & = & {\displaystyle \int_{a_1}^{+\infty}\frac{d}{dx}(f\circ \eta)(x)dx} \\
     & = & {\displaystyle \sum_{i=1}^{\infty}\int_{a_i}^{a_{i+1}}\frac{d}{dx}(f\circ \eta)(x)dx}\\
     & = & {\displaystyle \sum_{i=1}^{\infty}\int_{a_i}^{a_{i+1}}f_{x}(\eta(x))dx}. \\
 \end{array}
 $$
 This and $(f)$ imply that, for some sequence
 $x_n\to +\infty,$
 $$\lim_{n\to+\infty}f_{x_n}(\eta(x_n))=0.$$
 This contradiction, with the assumption that
 $\spec(Y)\cap (-\varepsilon,\varepsilon)=\emptyset,$
 proves the theorem.
\end{proof}

 \mm
 \section{Proof of Theorems 1.2 and 1.3}

 \mm
  To prove Theorem 1.2 we shall need the following.

\begin{lemma}[Lemma 6.2.11 of \cite{Ess}]\label{graded} Let $A=A_0
\oplus A_1 \oplus\cdots$ be a graded ring (A need not be
commutative). Let $a \in A_d$, for some $d \geq 1$. Then $1+a$ is
invertible in $A$ if and only if a is nilpotent.
\end{lemma}

\begin{proof}[\textbf{Proof of Theorem 1.2}]
We start as in the proof of Proposition 8.1.8 of \cite{Ess}. By the
Reduction Theorem (See \cite{BCW}, \cite{Dru}, \cite{Yag}) it
suffices to prove the Jacobian Conjecture for all $n \geq 2$ and all
polynomial maps $F:\mathbb{C}^n \rightarrow \mathbb{C}^n$ of the
form
$$F=(-X_1+H_1,\ldots,-X_n+H_n)$$
where $X_i:\mathbb{C}^n\rightarrow\mathbb{C}$ denotes the canonical
i-coordinate function, each $H_i$ is either zero or homogeneous of
degree 3 and $JH$ (with $H=(H_1,H_2,\ldots,H_n)$) is nilpotent.
Consider the polynomial map $\tilde{F}:\mathbb{R}^{2n} \rightarrow
\mathbb{R}^{2n}$ defined by
$$\tilde{F}=(ReF_1,ImF_1,\ldots,ReF_n,ImF_n).$$
So we have $\tilde{F}=-\tilde{X}+\tilde{H}$, where $\tilde{H}$ is
homogeneous of degree 3. Since $JH$ is nilpotent, $JF=-I+JH$ is
invertible by Lemma \ref{graded} and $\,\mbox{det}\,J\tilde{F}=
|\,\mbox{det}\, JF|^2=1$, whence $J\tilde{F}$ is invertible. So by
Lemma \ref{graded}, $J\tilde{H}$ is nilpotent and consequently
$\spec(\tilde{F})=\{-1\}$.

Now we proceed as in the proof of Proposition 8.3 of \cite{Jelonek}.
By \cite{Jel7} and \cite{Jel8}, we get that the set $S_{F}$ has
complex codimension 1, hence $S_{\tilde{F}}$ has  real codimension
2. Now $F$ is bijective if, and only if, $\tilde{F}$ is bijective,
Therefore if the assumption of this  theorem are satisfied, $F$ will
be bijective
\end{proof}

 \mm
 To prove Theorem 1.3 we shall need the following Jelonek results \cite{Jelonek}:

 \begin{theorem}\label{theo1Jelonek}
 if $Y:\mathbb{R}^n \rightarrow \mathbb{R}^n$ is a
 real polynomial mapping with nonzero Jacobian everywhere and
 ${\rm codim}(S_Y) \geq 3,$ then $Y$ is a bijection (and consequently $S_Y=\emptyset$).
 \end{theorem}

 \mm
 \begin{theorem}
 \label{theo2Jelonek}
 Let $Y:\mathbb{R}^n\rightarrow\mathbb{R}^m$ be a non-constant
 polynomial mapping. Then the set $S_Y$ is closed, semi-algebraic
 and for every non-empty connected component $S \subset S_Y$ we have
 $1\leq \dim (S) \leq n-1$. Moreover, for every point $q\in S_Y$ there
 exists a polynomial mapping $\phi: \mathbb{R}\rightarrow S_Y$ such
 that $\phi(\mathbb{R})$ is a semi-algebraic curve containing
 $\{q\}.$
 \end{theorem}

 \mm
 The proof of the following lemma is easy and will be omitted.

 \begin{lemma}\label{traslado} Let $Y:\mathbb{R}^n \rightarrow
 \mathbb{R}^n$ be a $C^1$-map such that \,$\spec(Y)\cap\{0\}=\emptyset$.
 Let $A:\mathbb{R}^n \rightarrow \mathbb{R}^n$ be a linear
 isomorphism. If $Z=A\circ Y \circ A^{-1}$  then $\spec(Y)=\spec(Z)$
 and $S_Z=A(S_Y)$.
 \end{lemma}

 \mm
 \begin{proposition}\label{a-main}
 Let $Y=(f,g,h):\mathbb{R}^3 \rightarrow \mathbb{R}^3$ be a polynomial map
 such that $\spec(Y) \cap (-\varepsilon,\varepsilon)=\emptyset$, for some
 $\varepsilon>0$. If ${\rm codim}(S_Y) \geq 2$ then $Y$ is a bijection.
 \end{proposition}

 \begin{proof}
 Suppose that $Y$ is not bijective. By Theorem \ref{theo1Jelonek}, we
 must have  $\dim(S_Y)=1$. Then by Theorem \ref{theo2Jelonek}, we
 obtain that

 \mm

 \begin{enumerate}
 \item[(a)] $Y(\mathbb{R}^3) \supset\mathbb{R}^3\setminus S_Y\,.$
 \end{enumerate}

 \mm
 Therefore, using again Theorem \ref{theo1Jelonek}, and
 Lemma \ref{traslado}, we way suppose that $S_Y$ contains a regular
 curve meeting transversally the plane $\{x=a\}$ at the point
 $p=(a,b,c)$. In this way.

 \mm
 \begin{enumerate}
 \item[(b)] the plane of $ \{x=a\}$
  contains a smooth embedded disc $D(a)$ such that $\{p\} = D(a) \cap S_Y$ and
 $C(a)\cap S_Y =\emptyset,$ where $C(a)$ is the boundary of $D(a).$
 \end{enumerate}

 \mm
 It is well known
 that there exists a positive integer $K$ such that

 \mm

 \begin{enumerate}
 \item[(c)] for all $q\in\mathbb{R}^3, \, \, \# Y^{-1}(q) \leq K.$
 \end{enumerate}

 \mm
 This implies that $Y^{-1}(C(a))$ is the union of finitely many embedded
 circles $C_1,C_2,\ldots,C_k$ contained in $f^{-1}(a)$. Each
 $Y_{|C_i}:C_i\rightarrow C(a)$ is a finite covering. As, by Theorem
 \ref{thm:fol}, each connected component of $f^{-1}(a)$ is a plane, we have that,
 for all $i=1,2,\ldots,k$, there exists a compact disc $D_i\subset
 f^{-1}(a)$ bounded by $C_i\,$.  It follows that, for all
 $i=1,2,\ldots,k$, $Y(D_i)=D(a)$. As $D(a)$ is simply connected, for all
 $i \in \{1,2,\ldots,k\}$, $Y_{|_{D_{i}}}:D_i\rightarrow D(a)$ is a
 diffeomorphism. Hence, if $q \in C$, $\# Y^{-1}(q)=k$.
 As $D(a)\cap S_{Y}=\{p\}$ and $\# Y^{-1}$ is locally constant, $\# Y^{-1}$
 must be identically equal to $k$ in $D(a)\setminus\{p\}$ and therefore $Y^{-1}(D(a)-\{p\})
 \subset\cup^{k}_{i=1}D_i$. As $Y$ is a local diffeomorphism, by using a limiting
 procedure,

 \mm
 \begin{enumerate}
 \item[(d)] for all $q \in D(a)$, $\# Y^{-1}(q)=k$ and so $Y^{-1}(D(a))=\cup^k_{i=1}D_i$.
 \end{enumerate}

\medskip
Notice that $D(a)$ can be taken of the form $D(a)=\{a\}\times D$,
where $D$ is a 2-disc of $\mathbb{R}^2$ centered at $(b,c)$; in this
way $C(a)=\{a\}\times\partial D$. We have that there exists
$\varepsilon > 0$ small that.

\mm
 \begin{enumerate}
 \item[(e)] if $s \in [a-\varepsilon,a+\varepsilon]$, $D(s)=\{s\}\times
 D$ and $C(s)=\{s\}\times\partial D$, then $(s,b,c)=D(s)\,\cap\,S_Y$
 and $C(s)\,\cap\,S_Y=\emptyset$.
 \end{enumerate}

\mm

Proceeding as above, we way find that for all $s \in
[a-\varepsilon,a+\varepsilon]$ there are $k$ embedded circles
$C_1(s),C_2(s),\ldots,C_k(s)$, with
$C_1(a)=C_1,C_2(a)=C_2,\ldots,C_k(a)=C_k$, contained in $f^{-1}(s)$
and such that $Y^{-1}(C(s))=\displaystyle{\cup^k_{i=1}}C(k(s))$.
Moreover each $C_i(s)$ depends continuously on $s$. Therefore,

\mm

 \begin{enumerate}
 \item[(f)] for all $s\in[a-\varepsilon,a+\varepsilon]$ and for all
 $i=1,2,\ldots,k$, there exists a compact disc $D_i(s) \subset
 f^{-1}(s)$ bounded by $C_i(s)$ such that $Y(D_i(s))=D(s)$ and
 $D_i(s)$ depends continuously on $s$.
 \end{enumerate}

\mm

Proceeding as in the proof of (d) we obtain that

 \begin{enumerate}
 \item[(g)] for all $s \in [a-\varepsilon,a+\varepsilon]$ and for all $q
 \in D(s)$, $\#Y^{-1}(q)=k$ and
 $Y^{-1}(D(s))=\displaystyle{\cup^k_{i=1}}D_i(s)$.
 \end{enumerate}

 As $[a-\varepsilon,a+\varepsilon]\times D$ is a compact neighborhood of
 $(a,b,c)$ and $Y^{-1}([a-\varepsilon,a+\varepsilon] \times D)$ is
 compact we obtain a contradiction with the assumption $p \in S_Y$.
 \end{proof}

 \mm
 The proof of the following  lemma can be found in \cite{agr} and
 \cite{fe-gu}. We include it here for sake of completeness.

 \begin{lemma}\label{main}
    Let $F:\mathbb{R}^n \rightarrow \mathbb{R}^n$ be a  differentiable map such
    that $\det(F^{\prime}(x))\neq 0$ for all $x$ in $\mathbb{R}^n.$
    Given $t\in\mathbb{R}$, let  $F_{t}\colon \mathbb{R}^n \rightarrow \mathbb{R}^n$
    denote the map
    $F_{t}(x)=F(x)-tx.$
    If there exists a sequence $\{t_m\}$ of real numbers
    converging to $0$ such that every map
    $F_{t_m}\colon \mathbb{R}^n \rightarrow \mathbb{R}^n$
    is injective, then $F$ is injective.
\end{lemma}

 \begin{proof}
    Choose $x_1,x_2\in \mathbb{R}^n$ such that
    $F(x_1)=y=F(x_2).$ We will prove $x_1=x_2.$
    By the Inverse Mapping Theorem,
    we may find neighborhoods $U_1, U_2, V$ of
    $x_1, x_2, y,$ respectively, such that, for $i=1,2,$
    $F|_{U_i} : U_i \to V$ is a diffeomorphism and
    $U_1\cap U_2 =\emptyset.$ If
    $m$ is large enough, then $F_{t_m} (U_1) \cap F_{t_m}(U_2)$
    will contain a neighborhood $W$ of $y$. In this way,
    for all $w\in W$, $\# (F_{t_m}^{-1}(w)) \ge 2$. This
    contradiction with the assumptions, proves the lemma.
 \end{proof}

 \mm
 \begin{remark}
    {\rm
    Even if $n=1$ and the maps $F_{t_m}$ in Lemma~\ref{main} are smooth
    diffeomorphisms, we cannot conclude that $F$ is a diffeomorphism.
    For instance, if $F:\mathbb{R}\to (0,1)$ is an orientation reversing
    diffeomorphism, then for every $t>0$, the map $F_t:\mathbb{R}\to\mathbb{R}$
    (defined by $F_t(x)=F(x)-tx$) will be an
    orientation reversing global diffeomorphism.
    }
\end{remark}

\mm
 \no
 \textbf{Theorem 1.3.}
 \textit{
 Let $Y=(f,g,h):\mathbb{R}^3 \rightarrow \mathbb{R}^3$ be a polynomial map
 such that $\spec(Y) \cap [0,\varepsilon)=\emptyset$, for some
 $\varepsilon>0$. If ${\rm codim}(S_Y) \geq 2$ then $Y$ is a bijection.
 }

 \begin{proof}

\no We claim that for each $0<t<\epsilon,$ the map $Y_t \colon
\mathbb{R}^3 \rightarrow \mathbb{R}^3,$ given by $Y_t(x)=Y(x)-tx,$
is injective.

In fact, as \, $D(Y_t)(x) = DY(x) - t I$, (where $I$ is the Identity
map), we obtain that if $0< a < \min \{t,\epsilon-t \} ,$ then
$\spec(Y_t)\cap (-a,a)=\emptyset.$ It follows immediately from
Lemma~\ref{main} and Proposition~\ref{a-main} that $Y$ is injective.
The conclusion of this theorem is obtained by using
Bia{\l}ynicki-Rosenlicht Theorem \cite{B-R}.
\end{proof}

 %%%%%%%%%%%%%%%%%%%%%%%%%%%%%%%%%%%%%%%%%%%%%%%%%%%%%%%%%%%%%%%%%%%%%%%%%%%%%%
%%%%%%%%%%%%%%%%%%%%%%%%%%%%%%%%%%%%%%%%%%%%%%%%%%%%%%%%%%%%%%%%%%%%%%%%%%%%%%%


\begin {thebibliography}{99}
\bibitem{alex}
    {\sc V. A. Alexandrov}, {\it Remarks on Efimov's Theorem about
    Differential Tests of Homeomorphism, \/} Rev. Roumanie Math. Pures
    Appl., {\bf 36} (1991), 3--4, pp. 101--105.


\bibitem{BCW}
    {\sc H. Bass, E.Connell and D. Wright}, {\it The Jacobian Conjecture:
    Reduction of Degree and Formal Expansion of the Inverse, \/}
    Bulletin of the American Mathematical Society, \textbf{7} (1982),
    287--330.


\bibitem{B-R} {\sc A. Bia{\l}ynicki-Birula and M. Rosenlicht}, {\it
Injective Morphisms of real algebraic varieties, \/} Proc. Amer.
Math. Soc., \textbf{13} (1962), 200--204.

\bibitem{Camacho}
    {\sc C. Camacho and Alcides Lins Neto}, {\it Geometric theory of foliations, \/}
    Birkh\"auser Boston Inc., (1985).

\bibitem{Cam}
    {\sc L.A. Campbell}, {\it Unipotent Jacobian matrices and
    univalent maps, \/} Contemp. Math. {\bf 264} (2000), 157--177.


%%\bibitem{Chamb-Meis}
%%    {\sc M. Chamberland and G. Meisters},
%%    {\it A mountain pass to the Jacobian Conjecture},
%%    Canadian Math. Bull. {\bf 41} (1998), 442--451.
%%\bibitem{CEGH}
%%    {\sc A. Cima, A. van den Essen, A. Gasull, E. Hubbers,  and F.
%%    Ma\~nosas,} {\it A Polynomial Counterexample to the Markus-Yamabe
%%    Conjecture, \/} Advances in Mathematics, {\bf 131} (1997),
%%    453--457.
%%\bibitem{CGM}
%%    {\sc A. Cima, A. Gasull}, and F. Ma\~nosas, {\it The Discrete
%%    Markus-Yamabe Problem. \/} Nonlinear Anal. {\bf 35} (1999),
%%   343--354.
%%

\bibitem{Cha}
    {\sc M. Chamberland}, {\it Characterizing two-dimensional maps
    whose Jacobians have constant eigenvalues.} Canad. Math. Bull. \textbf{46}
    (2003), no. 3, 323--331.

\bibitem{CGL}
  {\sc M. Cobo, C. Gutierrez, J. Llibre,} {\it On the injectivity of
   $C^1$ maps of the real plane } Canadian Journal of Mathematics
  {\bf 54} (2002),  No \textbf{6}, 1187--1201.

\bibitem{C-N}
  {\sc Chau and Nga} {\it A remark on Yu's theorem}, preprint,
  (1998).

\bibitem{Dru}
  {\sc L. Druskowski} {\it An effective approach to Keller's Jacobian conjecture},
  Math. Ann. \textbf{264} (1983), 303--313.

\bibitem{Ess}
{\sc A. van den Essen}, {\it Polynomial automorphisms and the
Jacobian conjecture, \/} Progress in Mathematics {\bf 190},
Birkhauser Verlag, Basel, (2000).

\bibitem{agr}
    {\sc A. Fernandes, C. Gutierrez and R. Rabanal,}
    { \it Global asymptotic stability for differentiable  vector fields of
 $ \mathbb{R}^2$.}  Journal of Differential Equations \textbf{206}, (2004), 470--482.


    \bibitem{fe-gu}
    {\sc A. Fernandes, C. Gutierrez and R. Rabanal,}
    { \it On local diffeomorphisms of $\mathbb{R}^n$ that are injective.}
    Qualitative Theory of Dynamical Systems \textbf{5}, (2004), 129--136, Article No. 63.

 \bibitem{godbillon}
    {\sc C. Godbillon},
    {\it Feuilletages: \'Etudes g\'eom\'etriques}
    Progress in Math. {\bf 14}, Birkh\"auser Verlag, (1991).

\bibitem{golu}
    {\sc O. Golubitsky and V. Guillemin},
    {\it Stable mappings and their singularities.}
    Grad. Texts in Math. {\bf 14}, Springer--Verlag, (1973).

\bibitem{Gu}
    {\sc C. Gutierrez,} {\it A Solution to the Bidimensional Global
    Asymptotic Stability Conjecture.} Ann. Inst. H. Poincar\'e.
    Analyse non Lineaire \textbf{12}, No.6 (1995) 627--671.

\bibitem{G-R} {\sc C. Gutierrez e R. Rabanal, \/}
{\it Injectivity of differentiable maps $\mathbb{R}^2\to
\mathbb{R}^2$ at infinity. \/} To appear in Bulletim of the Braz.
Math. Soc. (2006).

\bibitem{GV1} {\sc C. Gutierrez e Nguyen Van Chau. \/}
{\it  Properness and the Jacobian Conjecture in $\mathbb{R}^2.$ \/}
 Vietnam Journal of Mathematics. \textbf{31}, No. 4 (2003), 421--427.

\bibitem{GV2} {\sc C. Gutierrez e Nguyen Van Chau.  \/}
{\it  On Nonsingular Polynomial Maps of $\mathbb{R}^2.$ \/} To
appear in  Annales Polonici Mathematici. (2006).

\bibitem{GV3} {\sc C. Gutierrez e Nguyen Van Chau. \/}
{\it A remark on an eigenvalue condition for the global injectivity
of differentiable maps of $\mathbb{R}^2.$ \/ } Preprint.

\bibitem{Jelonek}  {\sc Z. Jelonek}, {\it Geometry of real polynomial mappings.\/}
     Math. Zeitschrift. \textbf{239} (2002), 321--333.

\bibitem{Jel7}
  {\sc Jelonek, Z} {\it The set of points at which a polynomial map is not
  proper.} Ann. Polon. Math. \textbf{58} (1993), 259--266.

\bibitem{Jel8}
  {\sc Jelonek, Z} {\it Testing sets for properness of polynomial mappings.}
  Math. Ann. \textbf{315} (1999), 1--35.

\bibitem{N-X}
    {\sc S. Nollet and F. Xavier,} {\it Global inversion via the
    Palais-Smale Condition,\/} Discret and Continuous Dynamical Systems
    Ser. A, {\bf 8} (2002), 17--28.

%%\bibitem{Ole}
%%    {\sc C. Olech,}{\it On the global stability of an autonomous
%%    system on the plane,\/} Cont. to Diff. Eq., {\bf 1}, (1963),
%%    389--400.

\bibitem{pal}
    {\sc C. F. B. Palmeira}, {\it Open manifolds foliated by planes}
    Ann. of Math. {\bf 107} (1978), 109--131.

\bibitem{CGPV} {\sc R. Peretz, Nguyen Van Chau, L. A. Campbell and C. Gutierrez,\/}
   {\it  Iterated Images and the Plane Jacobian Conjecture. \/}
   Discret and Continuous Dynamical Systems. \textbf{16}, No. 2 (2006), 455-461.

\bibitem{Pinchuk}
    {\sc S. Pinchuck}, {\it A counterexample to the strong Jacobian
    conjecture}, Math. Z. {\bf 217} (1994), 1--4.

\bibitem{Rosenberg}
    {\sc H. Rosenberg,} {\it Foliations by planes,} Topology,
    {\bf 7} (1968), 131--138.

\bibitem{S-X}
    {\sc B. Smyth, F. Xavier,} {\it Injectivity of local
    diffeomorphisms from nearly spectral conditions, \/} J. Diff.
    Equations, \textbf{130} (1996), 406--414.

\bibitem{Yag}
    {\sc A.V. Yagzhev, } {\it On Keller's problem, \/} Siberian Math.
    J., \textbf{21} (1980), 747--754.



\end {thebibliography}
\end{document}